\documentclass[12pt,authoryear]{elsarticle}

\usepackage{amssymb}
\usepackage{amsmath, amsthm}
\usepackage{graphicx}
\usepackage{color}
\usepackage{bm}
\usepackage{booktabs}
\usepackage{caption,subcaption}
\usepackage{multirow}
\usepackage{hyperref}
\usepackage{natbib}
\newtheorem{theorem}{Theorem}
\newtheorem{definition}{Definition}

\newtheorem{corollary}{Corollary}
\newtheorem{property}{Property}

\begin{document}

\begin{frontmatter}

\title{Expanding the Standard Diffusion Process to Specified Non-Gaussian Marginal Distributions}

\author[byu]{Robert Richardson\corref{cor1}}
\ead{richardson@stat.byu.edu}

\author[byu]{H. Dennis Tolley}
\ead{tolley@byu.edu}

\author[byu]{Kenneth Kuttler}
\ead{klkuttler@gmail.com}

\cortext[cor1]{Corresponding author}

\affiliation[byu]{
    organization={Department of Statistics, Brigham Young University},
    addressline={2152 WVB},
    city={Provo},
    postcode={84602},
    state={UT},
    country={USA}
}

\begin{abstract}
We develop a class of non-Gaussian translation processes that extend classical stochastic differential equations (SDEs) by prescribing arbitrary absolutely continuous marginal distributions. Our approach uses a copula-based transformation to flexibly model skewness, heavy tails, and other non-Gaussian features often observed in real data. We rigorously define the process, establish key probabilistic properties, and construct a corresponding diffusion model via stochastic calculus, including proofs of existence and uniqueness. A simplified approximation is introduced and analyzed, with error bounds derived from asymptotic expansions. Simulations demonstrate that both the full and simplified models recover target marginals with high accuracy. Examples using the Student’s t, asymmetric Laplace, and Exponentialized Generalized Beta of the Second Kind (EGB2) distributions illustrate the flexibility and tractability of the framework.
\end{abstract}


\end{frontmatter}

\section{Introduction}

Stochastic differential equations (SDEs) are widely used to model dynamic processes across various fields, including physical sciences, social sciences, and economics \citep{gardiner2009stochastic}. These equations have become indispensable tools for statisticians in data analysis, inference, and predictive modeling \citep{iacus2009simulation,kessler2012statistical,fuchs2013inference}. While traditional SDEs often rely on Brownian motion to drive the underlying dynamics, this paper introduces a novel approach that expands the standard diffusion model by replacing the Brownian motion process with a stochastic process that results in a user-specified, non-Gaussian marginal distribution.

We propose a framework where by using a copula-based transformation of Brownian motion we can achieve an error process with user-specific marginal distributions. The primary contribution of our work is the development of this prescriptive method while preserving the desirable properties of the diffusion framework.  This includes maintaining the decaying covariance structure characteristic of Brownian motion, but with greater flexibility. Unlike traditional methods, our approach bypasses the restrictive requirement that differences between successive observations conform to a model-ready distribution, offering a more adaptable framework for modeling complex phenomena.

Our model exhibits several key features that enhance its applicability across various domains:

\begin{enumerate}
    \item It encompasses the standard Gaussian process diffusion equation as a special case.
    \item The tail behavior, including the presence of fat tails, resembles processes with jumps, such as Lévy processes, providing a better fit for data with extreme events.
    \item The model captures longer memory effects similar to those in fractional Brownian motion, making it suitable for processes with long-term dependencies.
    \item Unlike fractional Brownian motion and 
     Lévy processes, the process does not require
     the error structure of the driving noise to be symmetric.
    \item The resulting process is not a martingale, except in the special case of Gaussian distributions, allowing for more flexible modeling of non-martingale behaviors.
    \item It does not assume independent increments, except in the Gaussian case, enabling the modeling of processes where this assumption does not hold.
    \item Despite these extensions, the process conditionally retains the Markov property, which is approximated by a random walk type process, ensuring that it remains tractable for analysis and simulation.
\end{enumerate}

Previous works have introduced specific error structures in place of Brownian motion. Two of the more popular variations on Brownian motion are Lévy processes and fractional Brownian motion. A Lévy process produces a diffusion process with jumps, meaning the resulting error process has a positive probability of a discontinuous jump \citep{bertoin1996Levy,schoutens2003Levy}. Fractional Brownian motion mimics the dependence structure of Brownian motion but with long term dependence and without the Markov property\citep{mandelbrot1968fractional,decreusefond1999stochastic}. Financial data especially has shown evidence of both long term dependence and jumps, making these alternatives popular for modeling such data. However, both approaches have limitations in their ability to directly model skewness in the error process. Lévy processes can generate asymmetry through jumps, but this skewness is not explicitly controlled or specified. Similarly, fractional Brownian motion retains a symmetric error structure due to its Gaussian nature, and does not inherently account for skewness.

Our approach not only mimics the long tail behavior of Lévy processes and the long memory behavior of fractional Brownian motion, but it directly addresses additional gaps by allowing the user to specify the marginal distribution of the process, including its skewness.  This capability is particularly important in fields like finance and risk management, where the asymmetry of returns or risks is a key feature of the data. While other models may capture some aspects of the phenomena we address, the flexibility to explicitly model skewness in the error process sets our approach apart, offering a more comprehensive tool for modeling real-world processes. To be clear, our proposed methodology does not hold Lévy processes and fractional Brownian motion as special cases, but it mimics their most attractive features.

There are a number of other alternatives that have been introduced. Many of these adjust the finite dimensional distribution of the resulting process, resulting in marginal distributions such as normal mixtures \citep{rombouts2014bayesian}, the inverse Gaussian distribution \citep{eriksson2009normal}, 
the inverse Gamma distribution \citep{milevsky1998asian}, the generalized Gamma distribution \citep{fabozzi2009estimating}, and others. 
Although these models produce marginal distributions in line with the specific distribution for which they are designed,
construction of each of these is specific to the marginal distribution of interest. Our approach has the distinct advantage that includes a unified framework that can produce any marginal absolutely continuous distribution. In contrast, the previous approaches required a unique framework for each model being used. The approach we will use may not necessarily hold all these other models as special cases by construction, but it will hold the marginal distributions of the processes described in these works as special cases. 

The following sections provide a detailed theoretical foundation for our method, followed by practical applications and case studies that demonstrate its utility in real-world scenarios. Section \ref{sec:theory} introduces the key elements of non-Gaussian translation processes and a number of key properties that show how it either conforms or differs from Brownian motion. Section \ref{sec:sde} takes a careful treatment of stochastic calculus using translation processes as the error process and evaluate approximation accuracy. Section \ref{sec:examples} builds on these results to show three specific distributional choices and confirm the approximations from section \ref{sec:sde} via simulation.  Section \ref{sec:conc} concludes with a discussion of the broader implications of our findings and approaches for future work.

\section{Non-Gaussian Translation Processes}
\label{sec:theory}
\subsection{Definition}

Let $(\Omega, \mathcal{F}, \mathbb{P})$ be a complete probability space. We consider a standard Brownian motion $\{B_t\}_{t \geq 0}$ defined on this space as a Gaussian process with $B_0 = 0$, $E[B_t] = 0$ and $Var[B_t] = t$ for $t\ge 0$. We denote by $\mathcal{F}_t$ the natural filtration generated by $B_t$. Define $F$ as an absolutely continuous (monotone increasing) distribution function corresponding to a member of a location-scale family with zero mean and unit variance. Define the sigma-algebra $\mathcal{G}$ to be the Borel sigma-algebra on $\mathbb{R}^n$.

We construct a general translation process $\{Z_t\}_{t \geq 0}$ through the following measurable transformation,
\begin{equation}
Z_t(\omega) = \sqrt{t} F^{-1} \left( \Phi \left( \frac{B_t(\omega)}{\sqrt{t}} \right) \right), \quad t > 0, \omega \in \Omega.
\label{trans}
\end{equation}
When $F$ is the Gaussian distribution function the process is
the Gaussian diffusion process.  Throughout this
paper the terms ``non-Gaussian" or ``translation process" indicates that $F$ is 
assumed to \emph{NOT} be a Gaussian distribution.

Practically this process extends the Brownian motion to a family of non-Gaussian processes by matching percentiles of the marginal distributions of the normal distribution to percentiles of the distribution function $F$.  This is a variant of the copula process methodology \citep{cherubini2004copula}
but used in a 
different way than is typically done. Copulas are typically constructed by placing a higher level correlation structure on top of 
transformed random variables, whereas the translation process places the transformation on a higher level than the correlation 
structure. The method also bears a resemblance to an Esscher and an Esscher-Girsanov transform for random variables
 \citep{goovaerts2008actuarial}. These transformations have been shown to capture a number of specific distributional 
 characteristics \citep{gerber1994option, lau2008option, monfort2012asset} but are not fully generalizable to any distribution.

For many of the properties we will discuss it is important that the distribution function $F$ corresponds to a standardized distribution in the first and second moments. The advantage of standardizing $F$ in this manner is that it allows flexibility in the marginal structure in the third and higher moments allowing for processes with skewness, heavy tails, etc.

\subsection{Properties}

We consider \( Z_t \) as a non-Gaussian translation of Brownian motion. The stochastic process \( \{ Z_t \}_{t \in T} \) is defined on a probability space \( (\Omega, \mathcal{G}, P^x) \), where \( Z_t: T \times \Omega \to \mathbb{R}^n \) is a measurable map. The probability measure \( P^x \) governs the finite-dimensional distributions
\[
\Pr[Z_{t_1} \in G_1, \ldots, Z_{t_k} \in G_k]
\]
for Borel sets \( G_i \subseteq \mathcal{G} \) in \( \mathbb{R}^n \). The measure \( P^x \) differentiates this process from Brownian motion, which is a Gaussian process with multivariate normal finite-dimensional distributions.

The translation process will have many of the same properties as Brownian motion and many that are different. We summarize a few of these properties and describe how they compare with Brownian motion
processes.  In all of these define $ Z_t $ as given in equation (\ref{trans}) to be a non-Gaussian translation process, where \( F \) is a non-Gaussian distribution function for a location-scale family of distributions. 
\begin{property}[Marginal Distribution]
The marginal distribution of \( Z_t \) at time \( t \) is given by \( F_{Z_t}(x) = F(x/\sqrt{t}) \).
\end{property}
\begin{property}[Moments]
If \( F \) is as defined above, with mean \( 0 \) and variance \( 1 \), then the first and second moments of \( Z_t \) at time \( t \) are \( \mathbb{E}[Z_t] = 0 \) and \( \text{Var}[Z_t] = t \).
\end{property}
Properties 1 and 2 are results of the construction of the process.

Let ${Y_k; k = 0,1,2,\dots}$ be a discrete process under the construction $Y_{k+\Delta t} = Y_{k} + \sqrt{\Delta t} \xi_k$ where $\xi_k$ is a random variable with support on $\{-1,1\}$ and having 0.5 probability of being either -1 or 1 with $\xi_k$ and $\xi_j$ being independent for $j \ne k$.   When $\Delta t$ is reasonably small, this process approximates Brownian motion well. An equivalent representation for the translation process is given in the following.
\begin{property}[Discrete Construction]
If \( F \) is \( C^1 \)-differentiable, then, \( Z_t \) can be constructed on discrete intervals using 
\begin{equation}
Z_{k+\Delta t} = Z_k + h(Y_k)\sqrt{\Delta t} \xi_k,
\label{Markov}
\end{equation}
where \( h(Y_k) = \frac{\phi (Y_k/\sqrt{k})}{f(F^{-1}(\Phi(Y_k/\sqrt{k})))} \) and $\xi_k$ is as defined above.  As \( \Delta t \to 0 \), we recover the continuous process defined by equation (\ref{trans}).
\end{property}
This essentially includes an adjustment factor that modifies the random walk. The function $h$ is  used in the proof for Property 3 given in Appendix A. Equation (\ref{Markov}) can be used to formulate the likelihood function of discrete time observations from this process.

\begin{property}[Martingale Property]
If \( F \) is \( C^\infty \)-smooth, then \( Z_t \) is NOT a martingale unless \( F \) is Gaussian.
\end{property}
\begin{proof}
Recall that
\[
Z_t = \sqrt{t}\, F^{-1}\Bigl[\Phi(B_t/\sqrt{t})\Bigr].
\]
Since \( Z_t \) is a one-to-one function of \( B_t \), we can express the conditional expectation as
\[
E[Z_t\mid Z_s] = E\Bigl[\sqrt{t}\, F^{-1}\Bigl(\Phi(B_t/\sqrt{t})\Bigr) \,\Bigm|\, B_s\Bigr].
\]
Define the function
\[
g(x) = \sqrt{t}\, F^{-1}\Bigl(\Phi(x/\sqrt{t})\Bigr).
\]
We now expand \( g(x) \) in a Taylor series about \( x = B_s \):
\[
g(B_t) = g(B_s) + g'(B_s)(B_t-B_s) + \frac{1}{2}g''(B_s)(B_t-B_s)^2 + \frac{1}{6}g'''(B_s)(B_t-B_s)^3 + \cdots.
\]
For convenience, define
\[
h(x) \equiv g'(x) = \frac{d}{dx}\Bigl[\sqrt{t}\, F^{-1}\Bigl(\Phi(x/\sqrt{t})\Bigr)\Bigr] 
= \frac{\phi(x/\sqrt{t})}{f\Bigl(F^{-1}(\Phi(x/\sqrt{t}))\Bigr)},
\]
where \(\phi\) denotes the standard normal density and \( f \) is the density corresponding to \( F \).

Taking the conditional expectation with respect to \( B_s \) and using the properties of Brownian motion (specifically, \( E[B_t-B_s \mid B_s]=0 \) and \( E[(B_t-B_s)^2 \mid B_s]=t-s \), with all odd powers vanishing), we obtain
\[
E[Z_t\mid Z_s] = g(B_s) + \frac{1}{2}g''(B_s)(t-s) + O\bigl((t-s)^{3/2}\bigr).
\]
Noting that \( g(B_s)=\sqrt{t}\, F^{-1}\Bigl[\Phi(B_s/\sqrt{t})\Bigr] \) and \( g''(B_s)=h'(B_s) \), we have
\[
E[Z_t\mid Z_s] = \sqrt{t}\, F^{-1}\Bigl[\Phi(B_s/\sqrt{t})\Bigr] + \frac{1}{2}h'(B_s)(t-s) + O\bigl((t-s)^{3/2}\bigr).
\]
If \( F \) is non-Gaussian, then in general \( h'(B_s) \neq 0 \) so that the term \( \frac{1}{2}h'(B_s)(t-s) \) introduces a drift. Consequently, \( Z_t \) is not a martingale. In contrast, if \( F \) is the Gaussian distribution then \( F^{-1}=\Phi^{-1} \) and \( h(x) \equiv 1 \), which implies \( h'(x)=0 \); in that case the process \( Z_t \) reduces to Brownian motion, a martingale.
\end{proof}

The assumption of differentiability required for the Taylor series expansion is met by the assumption that F is \( C^\infty \)-smooth. The Student's t, and the EGB2 distribution functions meet this assumption.  The asymmetric LaPlace is not \( C^\infty \)-smooth, in general.  However, as considered below, it does meet the conditions of differentiability.  Thus all the translation processes we use in this paper meet the requirements for this property.


\begin{property}[Non-stationary Increments]
For a translation process with \( s < t \), \( \text{Var}(Z_{t} - Z_s) = f(s,t) \neq f(t-s) \) if $F$ is not Gaussian, meaning that the translation process does not have stationary increments.
\end{property}
\begin{property}[Non-independent Increments]
For a translation process with \( s_1 < s_2  < t_1 < t_2 \), \( \text{Cov}(Z_{t_2} - Z_{t_1},Z_{s_2} - Z_{s_1}) \neq 0 \) if $F$ is not Gaussian, meaning that the translation process does not have independent increments.
\end{property}
Brownian motion famously has independent and stationary increments. It also has been shown that of processes with continuous paths, only Brownian motion has independent and stationary increments, so it is not surprising to see these two properties. The proof for both Properties 5 and 6 are shown in Appendix A. Not having stationary and independent increments, however, is an important property of fractional Brownian motion that is believed to make it better for modeling series of data with longer memories. 

Properties 4, 5 and 6 are not true if $F$ is a Gaussian distribution.

\begin{property}[Self-similarity]
Translation processes are self-similar, meaning they are invariant under scaling changes, i.e., \( \frac{1}{\sqrt{c}} Z_{ct} \overset{d}{=} Z_t \).
\end{property}
\begin{proof}
The proof of this is to show that distribution functions are the same: 
\begin{eqnarray*}
Pr\left(\frac{1}{\sqrt{c}} Z_{ct} < x\right) &=& Pr\left(\frac{1}{\sqrt{c}} \sqrt{ct} F^{-1} [\Phi(B_{ct}/\sqrt{ct})] < x \right) \\
&=& Pr\left( \sqrt{t} F^{-1} [\Phi(B_{ct}/\sqrt{ct})] < x \right) \\
&=& Pr\left(B_{ct}/\sqrt{c} < \sqrt{t}\Phi^{-1}[F(x/\sqrt{t})]\right)
\end{eqnarray*}
Because Brownian motion is self-similar, this is equivalent to $Pr\left(B_{t} < \sqrt{t}\Phi^{-1}[F(x/\sqrt{t})]\right)$ which can then be worked backwards through those same steps to determine the distribution function for $Z_t$. 
\end{proof}
Brownian motion and fractional Brownian motion are also self-similar. This is an important property for our application of the translation process as an error process in a diffusion SDE.

\subsection{The Cornish–Fisher Expansion }

An alternative approach to understanding the non-Gaussian corrections in the translation process is through the Cornish–Fisher expansion, which approximates the quantiles of a distribution based on a series expansion \citep{cornish1937moments}. Taking advantage of $F^{-1}(\Phi(B_t/\sqrt{t}))$ representing the $\Phi(B_t/\sqrt{t})$-th quantile of $F$, and the fact that we are forcing $F$ to have a mean of 0 and a standard deviation of 1, we can rewrite the transformation as
\begin{equation}
\sqrt{t}F^{-1}\Bigl(\Phi(B_t/\sqrt{t})\Bigr) = B_t + \sqrt{t}\delta(B_t/\sqrt{t})=Y_t.
\label{Corn delta}
\end{equation}
Set \(\kappa_3\) and \(\kappa_4\) to denote the standardized third and fourth cumulants (skewness and excess kurtosis) of $F$, respectively. The derivative of $Y_t$, \(\frac{\partial Y_t}{\partial \kappa_3}\), quantifies the sensitivity of the process to skewness.  In contrast, the derivative of $Y_t$, \(\frac{\partial Y_t}{\partial \kappa_4}\), measures the effect of excess kurtosis (heavy tails) on the process. In the Gaussian case, where \(\kappa_3=\kappa_4=0\), these derivatives vanish, and the process reduces to standard Brownian motion.

The deviation function \(\delta(x)\) given in
equation (\ref{Corn delta}) can be used to quantify the departure from Gaussian cases in the third and fourth cumulants. For a distribution with mean 0 and variance 1, the Cornish–Fisher expansion yields $\delta$ as
\[
\delta(x) = \frac{\kappa_3}{6}(x^2-1) + \frac{\kappa_4}{24}(x^3-3x) - \frac{\kappa_3^2}{36}(2x^3-5x) + \cdots.
\]

The function $\delta(x)$ allows us to see the effect of skewness and heavy tails on the process.  In many applications the higher-order cumulants are small, so retaining only the linear terms in \(\kappa_3\) and \(\kappa_4\) provides a reasonable approximation of the non-Gaussian effects.

For \(\kappa_3\), we obtain
\[
\frac{\partial Y_t}{\partial \kappa_3} = \sqrt{t}\,\frac{\partial \delta}{\partial \kappa_3}
=\sqrt{t}\cdot\frac{1}{6}\left[\left(\frac{B_t}{\sqrt{t}}\right)^2-1\right]
=\frac{1}{6\sqrt{t}}\,(B_t^2-t).
\]

Similarly, for \(\kappa_4\),
\[
\frac{\partial Y_t}{\partial \kappa_4} = \sqrt{t}\,\frac{\partial \delta}{\partial \kappa_4}
=\sqrt{t}\cdot\frac{1}{24}\left[\left(\frac{B_t}{\sqrt{t}}\right)^3-3\left(\frac{B_t}{\sqrt{t}}\right)\right]
=\frac{1}{24}\left(\frac{B_t^3}{t}-3B_t\right).
\]

A nonzero \(\kappa_3\) introduces asymmetry into the process, captured by the term \((B_t^2-t)\). Measures of excess kurtosis (heavy tails) are captured in the term \(\left(\frac{B_t^3}{t}-3B_t\right)\).

\section{Stochastic Calculus}
\label{sec:sde}

While the translation process is interesting in its own right, its main utility lies in serving as the driving noise process in stochastic differential equations (SDEs). Consider a general diffusion SDE driven by a translation process \(Z_t\):
\begin{equation}
dX_t = \alpha(X_t,t) \, dt + \sigma(X_t,t) \, dZ_t,
\label{sdeZ}
\end{equation}
where \(Z_t\) is the non-Gaussian translation of Brownian motion defined by
equation (\ref{trans})
and \(F\) is a distribution function from a location-scale family.  

\subsection{Existence of the Stochastic Integral}

\subsubsection{Preliminaries}

Brownian motion is defined as follows.

\begin{definition}
Let \( B_t \) be a stochastic process such that for all \( t_1 < t_2 < \cdots < t_m \), the increments \( \{ B_{t_i} - B_{t_{i-1}} \} \) are independent, and for any \( s < t \), the increment \( B_t - B_s \) is normally distributed with mean \( 0 \) and variance \( t - s \). Also, the sample paths \( t \mapsto B_t \) are almost surely H\"older continuous with every exponent \( \gamma \in (0, 1/2) \), , where \( \gamma \) denotes the H\"older continuity exponent. This means that for each \( \gamma < 1/2 \), there exists a constant \( C_\gamma \) such that
\[
|B_t - B_s| \leq C_\gamma |t - s|^\gamma
\]
for all \( s, t \) in a compact interval. This process is called a \emph{Wiener process} or \emph{Brownian motion}.
\end{definition}

One way to construct this process is via the expansion
\[
B_t(\omega) = \sum_{i=1}^\infty (\chi_{[0, t]}, g_i)_{L^2} \, \xi_i(\omega),
\]
where \( \{ \xi_i \}_{i=1}^\infty \) are independent standard normal random variables and \( \{ g_i \} \) is an orthonormal basis of \( L^2(0, \infty) \). Because Gaussian distributions have finite moments of all orders, for all \( \alpha > 0 \),
\[
\mathbb{E}\left( |B_t - B_s|^\alpha \right) \leq C_\alpha |t - s|^{\alpha/2} \mathbb{E}(|B_1|^\alpha).
\]
By the Kolmogorov--\v{C}entsov theorem \cite{KaratzasShreve1991,stroock2011probability}, this implies that \( t \mapsto B_t(\omega) \) is almost surely H\"older continuous with any exponent \( \gamma \in (0, 1/2) \).

However, Brownian motion is almost surely nowhere differentiable and not of bounded variation. In fact, the same is true for any non-constant continuous martingale. Consequently, Brownian motion cannot serve as an integrator in the classical Stieltjes sense for arbitrary continuous integrands. Instead, one must use the independence of the increments to define a stochastic integral of the form
\[
\int_0^t h(B_s, s) \, \sigma(X_s, s) \, dB_s.
\]
\index{Wiener process}This independence is also what makes $B_t$ a martingale and underlies the existence of the stochastic integral, which will also be a (local) martingale.

This integral then allows us to express stochastic differential equations (SDEs) in the form
\[
dX_t = \alpha(X_t, t) \, dt + h(B_t, t) \sigma(X_t, t) \, dB_t, \quad X(0) = X_0,
\]
where the last term is the stochastic integral. As in the case of the Stieltjes integral, its construction begins with approximating sums.

\subsubsection{Approximating the Stochastic Integral}

Let $C([0,T]; L^2(\Omega))$ denote the space of continuous functions from $[0,T]$ into $L^2(\Omega)$. Let $\Psi \in C([0,T]; L^2(\Omega))$ be adapted to the filtration $\{ \mathcal{F}_t \}$ generated by the Brownian motion $B$. We define a sequence of step function approximations:
\[
\Psi_n(t) = \sum_{i=0}^{n-1} \Psi(t_i) \, \chi_{(t_i, t_{i+1}]}(t),
\]
where $0 = t_0 < t_1 < \cdots < t_n = T$ is a partition of $[0,T]$ whose mesh goes to zero as $n \to \infty$. Define the approximate stochastic integral by
\[
\int_0^t \Psi_n \, dB := \sum_{i=0}^{n-1} \Psi(t_i) \left( B(t_{i+1} \wedge t) - B(t_i \wedge t) \right).
\]

Because of the independence of Brownian increments and Itô isometry, we have
\[
\mathbb{E}\left[ \left( \int_0^t \Psi_n \, dB \right)^2 \right] = \sum_{i=0}^{n-1} \mathbb{E}[\Psi(t_i)^2] (t_{i+1} \wedge t - t_i \wedge t).
\]
As $n \to \infty$, this converges to
\[
\int_0^t \mathbb{E}[\Psi(s)^2] \, ds.
\]
Hence $\int_0^t \Psi_n \, dB$ is a Cauchy sequence in $L^2(\Omega)$ and converges to a limit we denote
\[
\int_0^t \Psi \, dB.
\]

Moreover, using maximal inequalities and the Borel--Cantelli lemma, one can show that this limit has a version that is almost surely continuous in $t$. This defines the stochastic integral.

The stochastic integral is stable under perturbations of the integrand:
\[
\left\| \int_0^t \Psi \, dB - \int_0^t \Upsilon \, dB \right\|_{L^2(\Omega)}^2 \leq \int_0^T \mathbb{E}[(\Psi(s) - \Upsilon(s))^2] \, ds.
\]
It is also a martingale:
\[
\mathbb{E} \left[ \mathbf{1}_A \int_0^t \Psi \, dB \right] = \mathbb{E} \left[ \mathbf{1}_A \int_0^s \Psi \, dB \right] \quad \text{for } A \in \mathcal{F}_s, \, s < t.
\]

\subsubsection{Existence of a Solution}

\begin{theorem}\label{stochastic_integral}
Let $\Psi \in C([0,T]; L^2(\Omega))$ be adapted to the filtration $\{ \mathcal{F}_t \}$. Then the stochastic integral
\[
\int_0^t \Psi \, dB
\]
exists and is a continuous martingale. Moreover, it satisfies the Itô isometry:
\[
\mathbb{E}\left[ \left| \int_0^t \Psi \, dB \right|^2 \right] = \int_0^t \mathbb{E}[ \Psi(s)^2 ] \, ds.
\]
\end{theorem}

Although we assumed continuity of $\Psi$ in $t$, the integral can be defined for any $\Psi \in L^2([0,T]; L^2(\Omega))$ that is adapted. The stochastic integral remains a continuous martingale in this case.

Stochastic integrals naturally lead to stochastic differential equations. Consider the general form:
\[
dX_t = \alpha(X_t, t) \, dt + h(B_t, t) \, \sigma(X_t, t) \, dB_t, \quad X(0) = X_0,
\]
with $h$ continuous and possibly bounded.

\begin{theorem}\label{stochastic_DE}
Let \( \alpha \) and \( \sigma \) be Lipschitz in \( X \) and continuous in \( t \), and assume \( h \) is bounded and continuous. Then there exists a unique adapted process \( X \in C([0,T]; L^2(\Omega)) \) that solves
\[
dX_t = \alpha(X_t, t) \, dt + h(B_t, t) \sigma(X_t, t) \, dB_t, \quad X(0) = X_0.
\]
\end{theorem}

\begin{proof}
Define a mapping \( F \) by
\[
F(X)_t := X_0 + \int_0^t \alpha(X_s, s) \, ds + \int_0^t h(B_s, s) \sigma(X_s, s) \, dB_s.
\]
Endow \( C([0,T]; L^2(\Omega)) \) with the norm
\[
\| f \|_\nu := \sup_{t \in [0,T]} \| f(t) e^{-\nu t} \|,
\]
where \( \nu > 0 \) is a damping parameter introduced to control the growth of trajectories over time. 

Using the Lipschitz continuity of \( \alpha, \sigma \) and the boundedness of \( h \), one obtains
\[
\| FX - FY \|_\nu \leq \sqrt{\frac{C}{2\nu}} \| X - Y \|_\nu,
\]
for some constant \( C \) depending on \( T \), the Lipschitz constants, and \( \|h\|_\infty \). Choosing \( \nu \) large enough makes \( F \) a contraction. The Banach fixed point theorem guarantees a unique fixed point \( X \) which is adapted and continuous. This \( X \) solves the SDE almost surely.
\end{proof}

\begin{corollary}
In Theorem~\ref{stochastic_DE}, the boundedness assumption on $h$ can be relaxed. If $h$ is continuous, then the SDE still admits a unique adapted solution.
\end{corollary}

\begin{proof}
Let $\tau_n$ be the first exit time from the interval $[-n, n]$ for $|B_t|$. Define truncated versions of the SDE on $[0, \tau_n]$ with bounded $h_n := \chi_{[0,\tau_n]} h$, and let $X_n$ be the corresponding solutions.

By uniqueness, $X_n$ and $X_m$ agree on $[0, \tau_n]$ for $m > n$. Define $X(t) := \lim_{n \to \infty} X_n(t)$ pointwise. The limit exists, is adapted, and solves the original SDE almost surely.
\end{proof}

\subsection{Stochastic Integral Representation}

Having established the existence and uniqueness of the stochastic integral with a generalized integrand, we now consider how this integral can be approximated using discrete-time constructions and how it relates to the transformed Brownian motion defined earlier.

We aim to represent the stochastic integral \( \int_0^t \sigma(X_s,s) \, dZ_s \) using Riemann sums. However, caution is required: unless the integrator \( Z_t \) or the integrand \( \sigma(X_t, t) \) is of bounded variation, convergence of the Riemann sums in the pathwise (Riemann--Stieltjes) sense cannot be guaranteed. Since \( X_t \) involves a stochastic integral and is thus typically not of bounded variation, and since \( Z_t \) is a function of Brownian motion (which is almost surely nowhere differentiable), we cannot expect convergence in the classical sense.

That said, we can still assert convergence of the Riemann sums \emph{in probability}, which is sufficient for defining the stochastic integral. Specifically, we write
\begin{equation}
\int_0^t \sigma(X_s,s) \, dZ_s = \lim_{n \rightarrow \infty} \sum_{k=0}^{n-1} \sigma(X_{t_k}, t_k) \left( Z_{t_{k+1}} - Z_{t_k} \right),
\label{stoch_int}
\end{equation}
in the sense of convergence in probability. That is, for every \( \varepsilon > 0 \),
\[
\mathbb{P}\left( \left| \sum_{k=0}^{n-1} \sigma(X_{t_k}, t_k) (Z_{t_{k+1}} - Z_{t_k}) - \int_0^t \sigma(X_s,s) \, dZ_s \right| > \varepsilon \right) \to 0 \quad \text{as } n \to \infty.
\]

From Property 3, for small \(\Delta t = t_{k+1} - t_k\), we have
\begin{equation}
Z_{t_{k+1}} = Z_{t_k} + h(B_{t_k},t_k) \left( B_{t_{k+1}} - B_{t_k} \right),
\label{Markov_increment}
\end{equation}
where the function \(h\) is as defined above, namely
\begin{equation}
h(B_t,t) = \frac{\phi\left( \frac{B_t}{\sqrt{t}} \right)}{f\left( F^{-1} \left[ \Phi\left( \frac{B_t}{\sqrt{t}} \right) \right] \right)},
\label{h_function}
\end{equation}
with \(\phi\) and \(\Phi\) being the standard normal density and distribution functions, respectively, and \(f\) the density function corresponding to \(F\).

Substituting Equation \eqref{Markov_increment} into Equation \eqref{stoch_int}, we get
\begin{equation}
\int_0^t \sigma(X_s,s) \, dZ_s = \lim_{n \rightarrow \infty} \sum_{k=0}^{n-1} \sigma(X_{t_k}, t_k) \, h(B_{t_k},t_k) \left( B_{t_{k+1}} - B_{t_k} \right).
\label{stoch_int_B}
\end{equation}
This representation shows that the stochastic integral with respect to \(dZ_t\) can be expressed in terms of an integral with respect to \(dB_t\), for which the Itō integral is well-defined.

\subsection{Itō's Lemma for the Translation Process}
To further analyze \(Z_t\), we apply Itō's lemma to \(Z_t\) as a function of \(B_t\). Let us explicitly define
\[
g(y,t) = \sqrt{t} \, F^{-1} \Bigl[ \Phi\Bigl(\frac{y}{\sqrt{t}}\Bigr) \Bigr],
\]
where \(y\) represents the value of \(B_t\). The function \(g\) is chosen so that, by a standard change-of-variable argument, the exact process
\[
Z_t = g(B_t,t)
\]
has the desired marginal distribution (after appropriate scaling).

Applying Itō's lemma (see \cite{Oksendal2003}), we obtain
\[
dZ_t = \Bigl( \frac{\partial g}{\partial t} + \frac{1}{2} \frac{\partial^2 g}{\partial y^2} \Bigr) dt + \frac{\partial g}{\partial y}\, dB_t.
\]

Computing the derivatives explicitly, we have
\begin{align*}
\frac{\partial g}{\partial t} 
&= \frac{1}{2\sqrt{t}} \left[
    F^{-1}\Bigl( \Phi\Bigl(\frac{B_t}{\sqrt{t}}\Bigr) \Bigr)
    -  \, \frac{\frac{B_t}{\sqrt{t}}\, \phi\Bigl(\frac{B_t}{\sqrt{t}}\Bigr)}
    {f\Bigl(F^{-1}\Bigl(\Phi\Bigl(\frac{B_t}{\sqrt{t}}\Bigr)\Bigr)\Bigr)}
\right],\\[2mm]
\frac{\partial g}{\partial y} 
&= \frac{\phi\Bigl(\frac{B_t}{\sqrt{t}}\Bigr)}
{f\Bigl(F^{-1}\Bigl(\Phi\Bigl(\frac{B_t}{\sqrt{t}}\Bigr)\Bigr)\Bigr)},\\[2mm]
\frac{\partial^2 g}{\partial y^2} 
&= -\frac{1}{\sqrt{t}} \left[
     \, \frac{\frac{B_t}{\sqrt{t}}\, \phi\Bigl(\frac{B_t}{\sqrt{t}}\Bigr)}
    {f\Bigl(F^{-1}\Bigl(\Phi\Bigl(\frac{B_t}{\sqrt{t}}\Bigr)\Bigr)\Bigr)}
    + \frac{\phi^2\Bigl(\frac{B_t}{\sqrt{t}}\Bigr) \, f'\Bigl(F^{-1}\Bigl(\Phi\Bigl(\frac{B_t}{\sqrt{t}}\Bigr)\Bigr)\Bigr)}
    {f^3\Bigl(F^{-1}\Bigl(\Phi\Bigl(\frac{B_t}{\sqrt{t}}\Bigr)\Bigr)\Bigr)}
\right].
\end{align*}

Substituting these back into the expression for \(dZ_t\) and simplifying slightly, we have
\[
\begin{split}
dZ_t = \Biggl( &\frac{1}{2\sqrt{t}}\, F^{-1}\Bigl(\Phi\Bigl(\frac{B_t}{\sqrt{t}}\Bigr)\Bigr)
-\frac{1}{t}\, \frac{B_t\,\phi\Bigl(\frac{B_t}{\sqrt{t}}\Bigr)}
{f\Bigl(F^{-1}\Bigl(\Phi\Bigl(\frac{B_t}{\sqrt{t}}\Bigr)\Bigr)\Bigr)} \\
&\quad {} -\frac{1}{2\sqrt{t}}\,
\frac{\phi^2\Bigl(\frac{B_t}{\sqrt{t}}\Bigr)\, f'\Bigl(F^{-1}\Bigl(\Phi\Bigl(\frac{B_t}{\sqrt{t}}\Bigr)\Bigr)\Bigr)}
{f^3\Bigl(F^{-1}\Bigl(\Phi\Bigl(\frac{B_t}{\sqrt{t}}\Bigr)\Bigr)\Bigr)}
\Biggr)\, dt 
 + \frac{\phi\Bigl(\frac{B_t}{\sqrt{t}}\Bigr)}
{f\Bigl(F^{-1}\Bigl(\Phi\Bigl(\frac{B_t}{\sqrt{t}}\Bigr)\Bigr)\Bigr)}\, dB_t.
\end{split}
\]

\subsection{Rewriting the SDE}

Define the drift term as
\begin{equation}
\begin{split}
r(B_t,t) = \frac{1}{2\sqrt{t}}\, F^{-1}\Bigl(\Phi\Bigl(\frac{B_t}{\sqrt{t}}\Bigr)\Bigr)
-\frac{1}{t}\, \frac{B_t\,\phi\Bigl(\frac{B_t}{\sqrt{t}}\Bigr)}
{f\Bigl(F^{-1}\Bigl(\Phi\Bigl(\frac{B_t}{\sqrt{t}}\Bigr)\Bigr)\Bigr)}\\  -\frac{1}{2\sqrt{t}}\,
\frac{\phi^2\Bigl(\frac{B_t}{\sqrt{t}}\Bigr)\, f'\Bigl(F^{-1}\Bigl(\Phi\Bigl(\frac{B_t}{\sqrt{t}}\Bigr)\Bigr)\Bigr)}
{f^3\Bigl(F^{-1}\Bigl(\Phi\Bigl(\frac{B_t}{\sqrt{t}}\Bigr)\Bigr)\Bigr)}
\end{split}
\label{eq:remainder}
\end{equation}

Using Equation \eqref{Markov_increment}, we can rewrite the SDE in Equation \eqref{sdeZ} as
\begin{equation}
dX_t = \Bigl(\alpha(X_t,t) + r(B_t,t)\Bigr) \, dt + \sigma(X_t,t) \, h(B_t,t) \, dB_t.
\label{sdehZ_drift}
\end{equation}
This form allows us to leverage standard techniques for SDEs driven by Brownian motion, with the modification that the diffusion coefficient is now modulated by the function \(h(B_t,t)\), capturing the non-Gaussian characteristics of \(Z_t\).

Dropping the drift term yields a simplified stochastic differential equation (SDE)
\begin{equation}
\label{sdehZ}
dX_t = \alpha(X_t,t)\, dt + \sigma(X_t,t) \, h(B_t,t) \, dB_t.
\end{equation}
This simplified form offers several important advantages:

\begin{enumerate}
    \item \textbf{Asymptotic Negligibility:} 
    The drift correction involves increments with terms that are of order \( O(t^{-1/2}) \). As time increases, these terms decay and become negligible, so their exclusion introduces only a small approximation error in many practical settings.

    \item \textbf{Computational Simplicity:} 
    Without the drift term, there is no need to compute or approximate the derivative \( f' \) of the target density. This reduces the computational burden and avoids potential numerical inaccuracies that can arise when approximating derivatives.
    
    \item \textbf{Reduced Numerical Instability:} 
    The full drift term contains many small terms in both the numerator and denominator, which can lead to issues such as underflow (or overflow) during numerical computations. Although the function \( h(B_t,t) \) may also be sensitive to scaling, its numerical behavior is significantly more stable compared to the complete drift correction.
    
    \item \textbf{Easier Analytical Handling:} 
    Simplifying the SDE by omitting the drift correction makes the equation more tractable for analytical purposes. This can facilitate the derivation of further results or properties of the process, which might be obscured by the complexity of the full drift term.
    
\end{enumerate}

In summary, while the full drift term precisely adjusts the dynamics to capture the non-Gaussian characteristics of \( Z_t \), dropping it can lead to a more numerically robust, simpler, and often sufficiently accurate SDE for both simulation and theoretical analysis. In the next subsection we explore the implications of this approximation. Section \ref{sec:examples} includes a simulation study where the full SDE is compared against the approximated SDE.

\subsection{Error Analysis}

The remainder term has incremental terms according to \(r(B_t,t)\) which behave according to \(O(t^{-1/2})\). While dropping the drift term makes sense incrementally, we must explore the cumulative effects over time. 

We can write
\[
Z_t = Z_0 + \int_0^t h(B_s,s) \, dB_s + \int_0^t r(B_s,s) \, ds.
\]
Denote the approximate process by
\[
\widetilde{Z}_t = Z_0 + \int_0^t h(B_s,s) \, dB_s,
\]
where
\[
h(B_s,s)= \frac{\phi\Bigl(\frac{B_s}{\sqrt{s}}\Bigr)}{f\Bigl(F^{-1}\Bigl(\Phi\Bigl(\frac{B_s}{\sqrt{s}}\Bigr)\Bigr)\Bigr)}.
\]
Thus, the error (or remainder) incurred by approximating \(dZ_t\) solely by its stochastic term is
\[
E_t = Z_t - \widetilde{Z}_t = \int_0^t r(B_s,s) \, ds.
\]

Using the expressions for the derivatives, we observe that each term is individually \(O(t^{-1/2})\), noting that \(B_t / \sqrt{t}\) is \(O(1)\) with respect to \(t\). In other words, for sufficiently large \(t\) the drift term is dominated by an \(O(t^{-1/2})\) behavior, so that
\[
r(B_t,t)\,dt = O\Bigl(\frac{dt}{\sqrt{t}}\Bigr).
\]
In a discrete-time approximation with a small time step \(\Delta t\), the stochastic increment \(dB_t\) is \(O(\sqrt{\Delta t})\), while the drift term contributes an error of order \(O(\Delta t/\sqrt{t})\). Thus, as \(\Delta t/\sqrt{t} \to 0\), the error per step is negligible relative to the diffusion term. 

Over a fixed time horizon \(T\), the cumulative error (measured in an \(L^2\) sense) can be bounded by the Cauchy–Schwarz inequality:
\[
\mathbb{E}\left[ \Bigl| \int_0^T r(B_s,s)\,ds \Bigr|^2 \right] \le T \int_0^T \mathbb{E}\bigl[ r(B_s,s)^2 \bigr] ds.
\]
Assuming that \(\mathbb{E}[r(B_s,s)^2] = O(1/s)\), we obtain
\[
\mathbb{E}\left[ \Bigl| \int_0^T r(B_s,s)\,ds \Bigr|^2 \right] = O\Bigl(T \int_0^T \frac{ds}{s}\Bigr)
= O\Bigl(T \log T\Bigr).
\]
Thus, the root-mean-square accumulated error is
\[
\sqrt{\mathbb{E}\left[|E_T|^2\right]} = O\Bigl(\sqrt{T \log T}\Bigr).
\]
Since the typical fluctuations of the process (driven by the stochastic term) are \(O(\sqrt{T})\), the relative error is
\[
\frac{\sqrt{\mathbb{E}\left[|E_T|^2\right]}}{O(\sqrt{T})} = O\Bigl(\sqrt{\log T}\Bigr).
\]
This means that at its worst, the additional variation in the process added by the approximation is proportional to log time. 

It is important to note that the non-martingale nature of \(Z_t\) is exactly what introduces the drift term \(r(B_t,t)\). In the Gaussian case, where \(F^{-1}(\Phi(x)) = x\), the process is a martingale and \(r(B_t,t) \equiv 0\). For non-Gaussian \(F\), \(r(B_t,t)\) is nonzero and represents the deviation from martingale behavior. Also, we note that this is an upper bound. Empirical evidence has shown that dropping the drift term roughly keeps the shape of the expected marginal distribution \(f\), which will be demonstrated in the next section.

\section{Examples}
\label{sec:examples}

In this section, we propose three distributions that meet all the conditions required for the properties to hold, making a translation process useful as an error process in a diffusion model. These distributions are unique in that they allow for heavy tails, skewness, or both.

In the context of stochastic processes, Lévy processes provide examples of processes with heavy tails, meaning there is a non-negligible probability of observing large deviations. However, there are few examples of error processes that exhibit skewness. In the context of an error process, a skewed marginal distribution implies that errors in one direction occur more frequently, while errors in the opposite direction are less frequent but may have larger magnitudes. Typically, such asymmetry is captured in the drift component of a model rather than directly in the error process. By incorporating skewness directly into the error process through the translation process, we can model phenomena where the asymmetry is inherent in the random fluctuations rather than the deterministic trend.

\subsection{Illustrative Distributions}

Here we present three examples of candidate distributions that have parameters controlling skewness, tail width, or both. Recall that a condition for \( F \) in a translation process is that it must have mean 0 and variance 1. Each distribution discussed here will show how to constrain the parameters to meet that condition and how the parameters relate to different levels of tail width and skewness, which is an important feature of our expanded model.

\subsection*{Asymmetric Laplace Distribution}

The asymmetric Laplace distribution is a location-scale family distribution characterized by a skewness parameter \( \kappa > 0 \), a location parameter \( m \), and a scale parameter \( s \) \citep{yu2005three}. When \( \kappa = 1 \), the distribution is symmetric. Its probability density function (PDF) is given by:
\[
f(x; m, s, \kappa) = 
\begin{cases}
    \frac{\kappa}{s(1 + \kappa^2)} \exp\left(-\frac{\kappa (x - m)}{s}\right), & \text{if } x \geq m,\\
    \frac{\kappa}{s(1 + \kappa^2)} \exp\left(\frac{(x - m)}{s \kappa}\right), & \text{if } x < m.
\end{cases}
\]
To maintain a distribution with mean 0 and variance 1 for a given \( \kappa \), we set
\[
s = \left( \frac{\kappa^2}{1 + \kappa^4} \right)^{1/2}, \quad m = -s \frac{1 - \kappa^2}{\kappa}.
\]
By choosing or estimating a specific \( \kappa \), we adjust \( m \) and \( s \) accordingly to satisfy the required conditions. Here, \( \kappa \) controls skewness, with \( \kappa > 1 \) indicating positive skewness and \( \kappa < 1 \) indicating negative skewness.

\subsection*{Student's t-Distribution}

The Student's t-distribution is symmetric and has a parameter \( \nu \) (degrees of freedom) that controls the tail width. When \( \nu = 1 \), the distribution is a Cauchy distribution, and it approaches the normal distribution as \( \nu \rightarrow \infty \). The PDF of the t-distribution is given by:
\[
f(x; m, s, \nu) = \frac{\Gamma\left(\frac{\nu+1}{2}\right)}{\sqrt{\nu \pi} s \Gamma\left(\frac{\nu}{2}\right)} \left[ 1 + \frac{1}{\nu} \left( \frac{x - m}{s} \right)^2 \right]^{-\frac{\nu+1}{2}},
\]
where \( \Gamma(\cdot) \) is the gamma function. To use the t-distribution in a non-Gaussian transformation of Brownian motion, we include location and scale parameters \( m \) and \( s \), respectively. We set
\[
m = 0, \quad s = \left( \frac{\nu - 2}{\nu} \right)^{1/2}.
\]
Since the t-distribution has infinite variance when \( \nu \leq 2 \), we require \( \nu > 2 \) for the process to be well-defined. Here, \( \nu \) controls the heaviness of the tails, with smaller values indicating heavier tails.

\subsection*{Exponentialized Generalized Beta of the Second Kind}

The Exponentialized Generalized Beta of the Second Kind (EGB2) is a four-parameter distribution with location \( m \), scale \( s \), and shape parameters \( p \) and \( q \) \citep{mcdonald1995generalization}. Its PDF is given by:
\[
f(x; m, s, p, q) = \frac{|p|}{s \text{B}(p, q)} \left( \frac{\exp\left(\frac{x - m}{s}\right)}{1 + \exp\left(\frac{x - m}{s}\right)} \right)^{p} \left( 1 - \frac{\exp\left(\frac{x - m}{s}\right)}{1 + \exp\left(\frac{x - m}{s}\right)} \right)^{q-1},
\]
where \( \text{B}(p, q) \) is the beta function. The distribution is symmetric if \( p = q \), positively skewed if \( p > q \), and negatively skewed if \( p < q \). The tail width also varies with \( p \) and \( q \), with smaller values leading to heavier tails. To ensure the distribution has mean 0 and variance 1, we use
\[
s = \left( \Psi'(p) + \Psi'(q) \right)^{-1/2}, \quad m = \left( \Psi(q) - \Psi(p) \right) s,
\]
where \( \Psi \) is the digamma function, and \( \Psi' \) is its derivative. Here, the parameters \( p \) and \( q \) jointly control both skewness and tail behavior.

\subsection{Simulating Sample Paths}

Property 3 and Equations \eqref{sdehZ_drift} and \eqref{sdehZ} form the foundation for our modeling and inference approaches. In order to assess the importance of the drift correction in the random walk approximation, we simulate 10,000 sample paths using three different methods:
\begin{enumerate}
    \item \textbf{Standard Brownian motion}
    \item \textbf{Transformed Brownian motion} using Equation \eqref{trans}, which applies the non-Gaussian translation directly.
    \item \textbf{Random walk approximation} using both the drift-corrected version (Equation \eqref{sdehZ_drift}) and the uncorrected version (Equation \eqref{sdehZ}).
\end{enumerate}

Under standard Brownian motion, the marginal distribution is Gaussian, whereas the transformed processes are designed to follow the target non-Gaussian distributions. This simulation serves two purposes: first, to verify that our limiting approximations are appropriate; and second, to show that the random walk approximation—particularly when the drift correction is applied—closely captures the marginal distributions obtained by the direct transformation, a key requirement for our inference methods.

We simulate paths up to time \( t = 10 \) using a discretization with step size \( \Delta t = 0.01 \) (i.e., \( B_{t_k} = B_{t_{k-1}} + \sqrt{\Delta t}\, W_k \), where \( W_k \) are i.i.d. standard normal random variables). The simulation uses the following parameters for the distributions:
\begin{itemize}
\item \textbf{t-distribution (T)} with  (moderate tails) and  (very heavy tails), with location  and scale  ensuring unit variance.
\item \textbf{Asymmetric Laplace (AL)} distribution with  (moderate skew) and  (strong skew), with location  and scale  ensuring mean zero and unit variance.
\item \textbf{EGB2 distribution} with  (moderate skew and tail weight) and  (strong skew and heavier tails), with location  and scale  computed to ensure unit variance.
\end{itemize}

\begin{table}[htp]
    \centering
    \begin{tabular}{lccc}
        \toprule
        Distribution & Parameters & Location \( m \) & Scale \( s \) \\
        \midrule
        t-distribution & \( \nu = 10 \) & 0.000 & 0.894 \\
        t-distribution & \( \nu = 2.1 \) & 0.000 & 0.218 \\
        Asymmetric Laplace & \( \kappa = 1.5 \) & 0.508 & 0.609 \\
        Asymmetric Laplace & \( \kappa = 9 \) & 0.988 & 0.111 \\
        EGB2 & \( p = 0.95, q = 0.45 \) & -0.566 & 0.361 \\
        EGB2 & \( p = 4, q = 0.1 \) & (to be computed) & (to be computed) \\
        \bottomrule
    \end{tabular}
    \caption{Computed location \( m \) and scale \( s \) parameters ensuring mean zero and unit variance for each distribution. The EGB2 case for \( p = 4, q = 0.1 \) is pending computation.}
    \label{tab:parameters}
\end{table}

\begin{figure}[htp]
  \centering
  \begin{subfigure}[t]{0.45\textwidth}
    \centering
    \includegraphics[width=\textwidth]{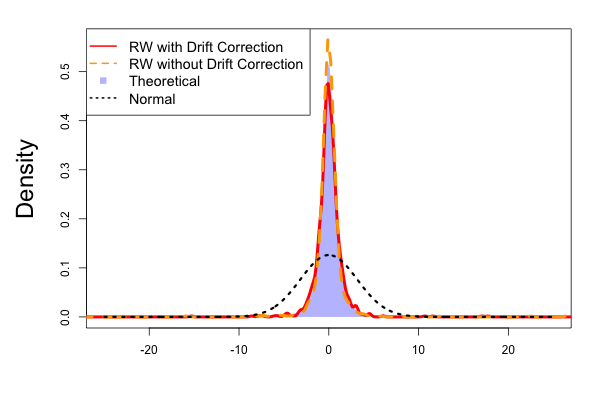}
    \caption[t-distribution, $\nu=2.1$]{T distribution with thin tails}
  \end{subfigure}
  \hfill
  \begin{subfigure}[t]{0.45\textwidth}
    \centering
    \includegraphics[width=\textwidth]{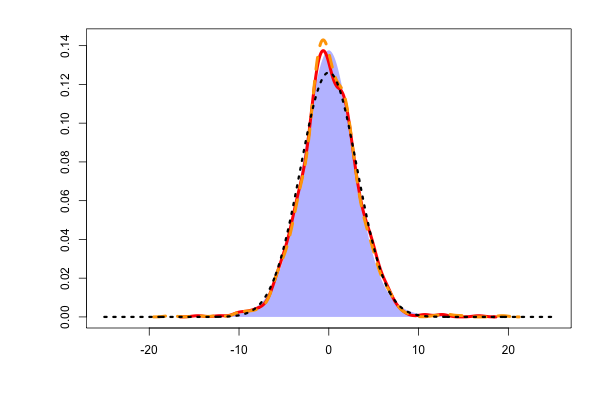}
    \caption[t-distribution, $\nu=10$]{T distribution with thick tails}
  \end{subfigure}

  \medskip

  \begin{subfigure}[t]{0.45\textwidth}
    \centering
    \includegraphics[width=\textwidth]{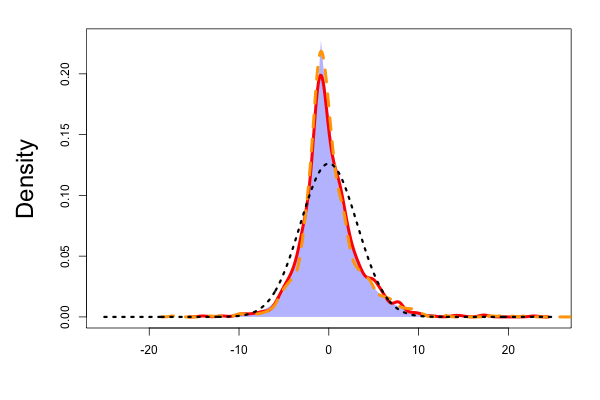}
    \caption[Asymmetric Laplace, $\kappa=1.5$]{Asymmetric Lalplace with a light skewness.}
  \end{subfigure}
  \hfill
  \begin{subfigure}[t]{0.45\textwidth}
    \centering
    \includegraphics[width=\textwidth]{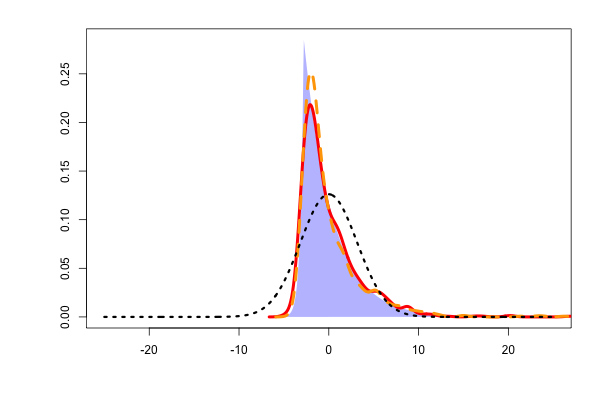}
    \caption[Asymmetric Laplace, $\kappa=9$]{Asymmetric Laplace with a strong skewness.}
  \end{subfigure}

  \medskip

  \begin{subfigure}[t]{0.45\textwidth}
    \centering
    \includegraphics[width=\textwidth]{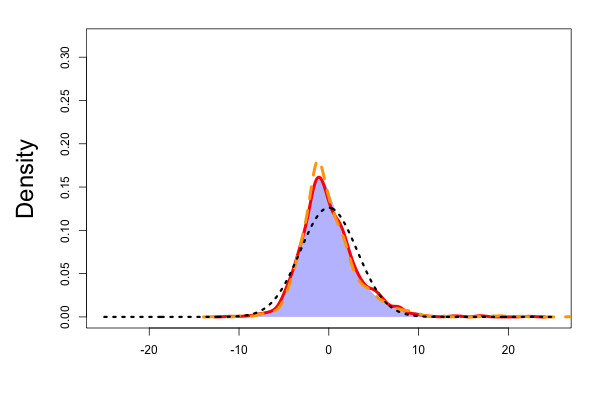}
    \caption[EGB2, $p=0.95,q=0.45$]{Log EGB2 with a weak skewness.}
  \end{subfigure}
  \hfill
  \begin{subfigure}[t]{0.45\textwidth}
    \centering
    \includegraphics[width=\textwidth]{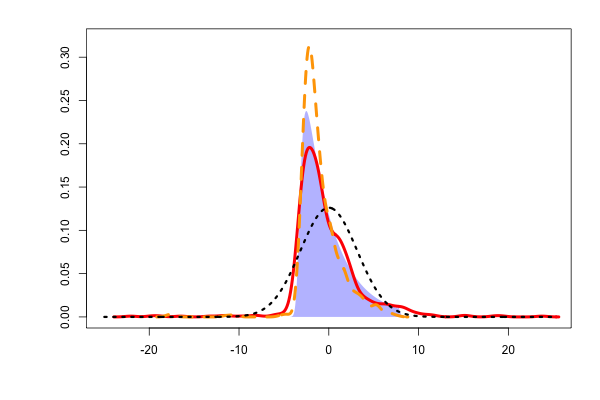}
    \caption[EGB2, $p=4,q=0.1$]{Log EGB2 with a strong skewness}
  \end{subfigure}

  \caption{Simulated marginal density functions at $t=10$ for 10,000 sample paths under different distributional assumptions. Each plot shows the theoretical density (shaded blue area), the drift-corrected approximation (red line), the approximation without drift correction (orange line), and a standard normal distribution (black dashed line) for reference.}
  \label{fig:combined}
\end{figure}

Figure \ref{fig:combined} displays the simulated marginal density functions at \( t = 10 \). The  random walk from Equation \eqref{sdehZ_drift} is shown on the same plot as the approximation from Equation \eqref{sdehZ} is shown. In each figure, the simulated density using the drift (solid red line) and without the drift (dashed orange line) is compared against the theoretical target (shaded blue area).   Also included is the corresponding normal distribution if Brownian Motion was used as the error process instead (black dashed line).

The red line should match the theoretical distribution exactly, which can be seen in Figure \ref{fig:combined}. It is able to match the thicker tails of the $t$ distribution, the skewness of the asymmetric Laplace distribution, and the skewness and thicker tails of the EGB2 distribution, even as we vary the parameters which modify these characteristics. The approximation shown by the orange line is not quite as exact, but is very close. The approximation also very clearly matches the tail width and skewness of the target distribution and follows the shape as the parameters change. 

We make the following conclusions from this demonstration:
\begin{enumerate}
    \item We are able to recover a specific marginal structure as the resulting sample path of a diffusion process
    \item Using the simpler approximation without the drift term, while not as precise, also matches the general marginal shape
\end{enumerate}

\section{Conclusion}
\label{sec:conc}

In this paper, we have introduced non‐Gaussian translation processes as a novel and versatile framework for extending stochastic differential equation models. This approach offers a straightforward yet creative means to incorporate user‐specified, non‐Gaussian marginal distributions into SDEs while preserving the desirable correlation structure driven by an underlying Brownian motion. A key advantage of our method is its flexibility: the marginal distribution can be tailored to match empirical data without sacrificing the tractability of the diffusion framework. Moreover, we have demonstrated that the translation process can be discretized into a random walk with an adjusted step size, thereby bridging the gap between traditional random walk models used in financial applications and continuous-time stochastic processes.

Through a basic analysis using Ito's lemma, we have shown that the resulting SDE representation of the translated process approximates the corresponding random walk dynamics for large 
$t$, providing both a mathematical foundation and practical guidance for implementation. Our simulation studies further illustrate that the method accurately reproduces the target non-Gaussian behavior while maintaining the underlying dependence structure. In addition to this, we have provided motivation for a simpler representation, which is also motivated by simulation results. 

There remain several promising directions for future work. One immediate avenue is the development of robust inference procedures for these processes, particularly for non-linear SDEs such as the Ornstein–Uhlenbeck or CIR models. Once these inferential procedures are fully developed, we can do more direct comparisons to models using Levy processes and fractional Brownian motion. Explicitly modeling and accounting for non-Gaussian marginal shapes may yield significant improvements in a wide range of applications, including finance, economics, and longitudinal data analysis. Furthermore, investigating alternative discretization schemes and higher-order corrections could enhance the numerical accuracy and stability of the proposed approach. 

Overall, non‐Gaussian translation processes represent a powerful tool for modeling complex stochastic phenomena. We believe that the flexibility and tractability of this framework will stimulate further research and find broad application in both theoretical and applied settings. 

\begin{appendix}

\section{Proofs of Section 2 Properties}
\label{rwproof}

\subsubsection*{Proof of Property 3}

Suppose that the random process has marginal distributions of the form $F_{\mu_t,\sigma_t}$ where $\mu_t$ is a location parameter and  $\sigma_t$ is a scale parameter. Other parameters may be present as well in $F$, such as degrees of freedom in a $t$ distribution.   

When a random walk is implemented with step sizes of $\Delta t$, then the distribution of the process at the $k$-th step is $F_{\mu_{k \Delta t},\sigma_{k \Delta t}}$. A translational non-Gaussian process can be expressed at step $k$ as 
\begin{equation}
Z_{k}=Z_{(k-1)}+F_{\mu _{k\Delta t},\sigma _{_{k\Delta t}}}^{-1}\left( \Phi
\left( Y_k\right) \right) -F_{\mu _{_{((k-1))\Delta t}},\sigma _{((k-1))\Delta
t}}^{-1}\left( \Phi \left( Y_{(k-1)}\right) \right)  
\label{general eq}
\end{equation}
where $Y_k$ represents a scaled Brownian motion, meaning that $Y_k = \frac{1}{\sqrt{k}} \sum_{i=1}^k \xi_i$ where each $\xi_i$ is an independent binary random variable which, for every $i$, is either 1 or -1, with probability .5 of each. 

Equation (\eqref{general eq}) can be rewritten using the property that $F_{\mu,\sigma}^{-1} = \sigma F_{0,1}^{-1} + \mu$, to give
\[Z_{k}=Z_{(k-1)}+\sigma _{k\Delta t}F_{0,1}^{-1}\left( \Phi \left( Y_{k}\right)
\right) -\sigma _{((k-1))\Delta t}F_{0,1}^{-1}\left( \Phi \left(
Y_{(k-1)}\right) \right) +\mu _{k\Delta t} - \mu _{((k-1))\Delta t}.
\]
The central piece of this equation can be rewritten as \small
\begin{eqnarray*}
\sigma _{k\Delta t}F_{0,1}^{-1}\left( \Phi \left( Y_{k}\right)
\right) &-&\sigma _{((k-1))\Delta t}F_{0,1}^{-1}\left( \Phi \left(
Y_{(k-1)}\right) \right) = \\ && \sigma_{k\Delta t}\left(F_{0,1}^{-1}\left( \Phi \left( Y_{k}\right)\right) - F_{0,1}^{-1}\left( \Phi \left(
Y_{(k-1)}\right)\right) +\left(1-\frac{\sigma_{((k-1))\Delta t}}{\sigma_{k\Delta t}}\right) F_{0,1}^{-1}\left( \Phi \left(
Y_{(k-1)}\right)\right) \right) 
\end{eqnarray*}

If we require the variance process to follow $\lim_{k \rightarrow \infty} \frac{\sigma_{((k-1))\Delta t}}{\sigma_{k\Delta t}} = 1$, which is a typical characteristic of Brownian motion, then the right side of the equation (\eqref{general eq}) can be approximated by
\[Z_{k}=Z_{(k-1)}+\ \sigma_{k\Delta t}\left(F_{0,1}^{-1}\left( \Phi \left( Y_{k}\right)\right) - F_{0,1}^{-1}\left( \Phi \left(
Y_{(k-1)}\right)\right)\right) +\mu _{k\Delta t} - \mu _{((k-1))\Delta t}.
\]

The following two identities will be used to express the translational non-Gaussian process as a random walk:
\begin{eqnarray}
\lim_{n\rightarrow \infty }\sqrt{n-1}\left( f(x)-f\left( \frac{x\sqrt{n}-1}{%
\sqrt{n-1}}\right) \right) &=& \frac{d}{dx}(f(x)) \label{id1} \label{id1} \\
\lim_{n\rightarrow \infty }\sqrt{n-1}\left( f(x)-f\left( \frac{x\sqrt{n}+1}{%
\sqrt{n-1}}\right) \right) &=& -\frac{d}{dx}(f(x)) \label{id2}
\end{eqnarray}

The scaled Brownian motion process can be written recursively and solved backwards when $\xi_k$ = 1 as $Y_{(k-1)} = \frac{Y_k \sqrt{k} - 1}{\sqrt{(k-1)}}$ and when $\xi_k = -1$ as $Y_{(k-1)} = \frac{Y_k \sqrt{k} + 1}{\sqrt{(k-1)}}$. Using this recursion with equations (\eqref{id1}) and (\eqref{id2}), the translational non-Gaussian process can be approximated as 
\begin{equation}
Z_{k} = Z_{(k-1)} +\mu _{k\Delta t} - \mu _{((k-1))\Delta t} + u(\sigma_{k\Delta t}, Y_{k}) \xi_k
\label{rwalmost}
\end{equation}
where \[u(\sigma_{k\Delta t}, Y_{k}) = \frac{1}{\sqrt{k}}\sigma _{k\Delta t} \frac{d}{dx}\left(
F_{0,1}^{-1}\left( \Phi \left( x\right) \right) \right) |_{x=Y_{k}}\]

This is equivalent to a standard random walk used to derive Brownian motion but with a correction factor $u(\sigma_{k\Delta t}, Y_{k})$. Note that if $\sigma _{k\Delta t} =  \sigma \sqrt{k \Delta t} $, then correction factor reduces to $u(\sigma_{k\Delta t}, Y_{k}) = \sigma  \sqrt{\Delta t} \frac{d}{dx}\left(
F_{0,1}^{-1}\left( \Phi \left( x\right) \right) \right) |_{x=Y_{k}}$. The derivative $\frac{d}{dx}\left(
F_{0,1}^{-1}\left( \Phi \left( x\right) \right) \right) |_{x=Y_{k}}$ is equal to $\frac{\phi(Y_k)}{f(F^{-1}_{0,1}(\Phi(Y_k)))}$ where $\phi$ is the standard normal density function and $f$ is the density function corresponding to the distribution function $F$. In this case, equation (\eqref{rwalmost}) simplifies to 
\begin{equation}
Z_{k} = Z_{(k-1)} +\mu _{k\Delta t} - \mu _{((k-1))\Delta t} + \sigma  \sqrt{\Delta t}\frac{\phi(Y_k)}{f(F^{-1}_{0,1}(\Phi(Y_k)))} \xi_k
\end{equation}

\subsubsection*{Proof of Property 5}

\begin{proof}
The difference $Z_t - Z_s$ has a mean of 0, so 
\begin{eqnarray*}
Var(Z_t - Z_s) &=& E[(Z_t - Z_s)^2] \\
&=& E[Z_t^2] + E[Z_s^2] - 2 E[Z_tZ_s] \\
&=& t + s - 2 E[Z_tZ_s]
\end{eqnarray*}
The last step comes from the construction of the translation process. Using the law of total expectation
The difference $Z_t - Z_s$ has a mean of 0, so 
\begin{eqnarray*}
E[Z_tZ_s] &=& E_{Z_s}\left[E_{Z_t}(Z_t Z_s | Z_s)\right] \\
&=&  E_{Z_s}\left[ Z_s E_{Z_t}(Z_t | Z_s)\right] \\
\end{eqnarray*}
The proof for Property 3 already showed that $E(Z_t|Z_s) = \sqrt{t} F^{-1} [\Phi(B_s/\sqrt{t})] + h'(B_s) (t-s)$. 
\begin{eqnarray*}
E[Z_tZ_s] &=&  E_{Z_s}\left[ Z_s \sqrt{t} F^{-1} [\Phi(B_s/\sqrt{t})] \right]  + E_{Z_s}\left[ Z_s h'(B_s) (t-s) \right]\\
\end{eqnarray*}
To provide a counterexample, let $B_s = 0$, which removes the second term due to  $h'(0) = 0$. 
\begin{eqnarray*}
E[Z_tZ_s]  &=& E_{Z_s}\left[ Z_s \sqrt{t} F^{-1} [\Phi(B_s/\sqrt{t})] \right]  \\
&=& \sqrt{st} F^{-1} (.5)\\
\end{eqnarray*}
Then for this case, $Var(Z_t - Z_s) = t + s - \sqrt{st} F^{-1} (.5) \neq f(t-s)$. 
\end{proof}

\subsubsection*{Proof of Property 6}
Because the processes all have mean 0, 
\begin{eqnarray*}
Cov(Z_{t_2} - Z_{t_1},Z_{s_2} - Z_{s_1}) &=& E(Z_{t_2} Z_{s_2}) + E(Z_{t_1} Z_{s_1})  - E(Z_{t_2} Z_{s_1}) - E(Z_{t_1} Z_{s_2})
\end{eqnarray*}
From the proof for Property 5 we see a case where $E(Z_t Z_s) = \sqrt{st} F^{-1}(.5)$. We will use this case again as a counterexample. Using this example the covariance is $F^{-1}(.5)(\sqrt{s_2t_2} + \sqrt{s_1t_1} - \sqrt{s_1t_2} - \sqrt{s_2t_1}) \neq 0$.

\end{appendix}

\bibliographystyle{elsarticle-harv}
\bibliography{transbib}

\begin{thebibliography}{23}
\expandafter\ifx\csname natexlab\endcsname\relax\def\natexlab#1{#1}\fi
\providecommand{\url}[1]{\texttt{#1}}
\providecommand{\href}[2]{#2}
\providecommand{\path}[1]{#1}
\providecommand{\DOIprefix}{doi:}
\providecommand{\ArXivprefix}{arXiv:}
\providecommand{\URLprefix}{URL: }
\providecommand{\Pubmedprefix}{pmid:}
\providecommand{\doi}[1]{\href{http://dx.doi.org/#1}{\path{#1}}}
\providecommand{\Pubmed}[1]{\href{pmid:#1}{\path{#1}}}
\providecommand{\bibinfo}[2]{#2}
\ifx\xfnm\relax \def\xfnm[#1]{\unskip,\space#1}\fi
\bibitem[{Bertoin(1996)}]{bertoin1996Levy}
\bibinfo{author}{Bertoin, J.}, \bibinfo{year}{1996}.
\newblock \bibinfo{title}{L{\'e}vy processes}. volume \bibinfo{volume}{121}.
\newblock \bibinfo{publisher}{Cambridge university press Cambridge}.
\bibitem[{Cherubini et~al.(2004)Cherubini, Luciano and Vecchiato}]{cherubini2004copula}
\bibinfo{author}{Cherubini, U.}, \bibinfo{author}{Luciano, E.}, \bibinfo{author}{Vecchiato, W.}, \bibinfo{year}{2004}.
\newblock \bibinfo{title}{Copula methods in finance}.
\newblock \bibinfo{publisher}{John Wiley \& Sons}.
\bibitem[{Cornish and Fisher(1937)}]{cornish1937moments}
\bibinfo{author}{Cornish, E.A.}, \bibinfo{author}{Fisher, R.A.}, \bibinfo{year}{1937}.
\newblock \bibinfo{title}{Moments and percentile points of distributions}.
\newblock \bibinfo{journal}{Supplement to the Journal of the Royal Statistical Society} \bibinfo{volume}{4}, \bibinfo{pages}{307--320}.
\bibitem[{Decreusefond and {\"U}st{\"u}nel(1999)}]{decreusefond1999stochastic}
\bibinfo{author}{Decreusefond, L.}, \bibinfo{author}{{\"U}st{\"u}nel, A.S.}, \bibinfo{year}{1999}.
\newblock \bibinfo{title}{Stochastic analysis of the fractional brownian motion}.
\newblock \bibinfo{journal}{Potential analysis} \bibinfo{volume}{10}, \bibinfo{pages}{177--214}.
\bibitem[{Eriksson et~al.(2009)Eriksson, Ghysels and Wang}]{eriksson2009normal}
\bibinfo{author}{Eriksson, A.}, \bibinfo{author}{Ghysels, E.}, \bibinfo{author}{Wang, F.}, \bibinfo{year}{2009}.
\newblock \bibinfo{title}{The normal inverse gaussian distribution and the pricing of derivatives}.
\newblock \bibinfo{journal}{The Journal of Derivatives} \bibinfo{volume}{16}, \bibinfo{pages}{23--37}.
\bibitem[{Fabozzi et~al.(2009)Fabozzi, Tunaru and Albota}]{fabozzi2009estimating}
\bibinfo{author}{Fabozzi, F.J.}, \bibinfo{author}{Tunaru, R.}, \bibinfo{author}{Albota, G.}, \bibinfo{year}{2009}.
\newblock \bibinfo{title}{Estimating risk-neutral density with parametric models in interest rate markets}.
\newblock \bibinfo{journal}{Quantitative Finance} \bibinfo{volume}{9}, \bibinfo{pages}{55--70}.
\bibitem[{Fuchs(2013)}]{fuchs2013inference}
\bibinfo{author}{Fuchs, C.}, \bibinfo{year}{2013}.
\newblock \bibinfo{title}{Inference for Diffusion Processes: With Applications in Life Sciences}.
\newblock \bibinfo{publisher}{Springer Science \& Business Media}.
\bibitem[{Gardiner(2009)}]{gardiner2009stochastic}
\bibinfo{author}{Gardiner, C.}, \bibinfo{year}{2009}.
\newblock \bibinfo{title}{Stochastic methods}. volume~\bibinfo{volume}{4}.
\newblock \bibinfo{publisher}{springer Berlin}.
\bibitem[{Gerber et~al.(1994)Gerber, Shiu et~al.}]{gerber1994option}
\bibinfo{author}{Gerber, H.U.}, \bibinfo{author}{Shiu, E.S.}, et~al., \bibinfo{year}{1994}.
\newblock \bibinfo{title}{Option pricing by esscher transforms}.
\newblock \bibinfo{journal}{Transactions of the Society of Actuaries} \bibinfo{volume}{46}, \bibinfo{pages}{140}.
\bibitem[{Goovaerts and Laeven(2008)}]{goovaerts2008actuarial}
\bibinfo{author}{Goovaerts, M.J.}, \bibinfo{author}{Laeven, R.J.}, \bibinfo{year}{2008}.
\newblock \bibinfo{title}{Actuarial risk measures for financial derivative pricing}.
\newblock \bibinfo{journal}{Insurance: Mathematics and Economics} \bibinfo{volume}{42}, \bibinfo{pages}{540--547}.
\bibitem[{Iacus(2009)}]{iacus2009simulation}
\bibinfo{author}{Iacus, S.M.}, \bibinfo{year}{2009}.
\newblock \bibinfo{title}{Simulation and inference for stochastic differential equations: with R examples}.
\newblock \bibinfo{publisher}{Springer Science \& Business Media}.
\bibitem[{Karatzas and Shreve(1991)}]{KaratzasShreve1991}
\bibinfo{author}{Karatzas, I.}, \bibinfo{author}{Shreve, S.E.}, \bibinfo{year}{1991}.
\newblock \bibinfo{title}{Brownian Motion and Stochastic Calculus}. volume \bibinfo{volume}{113} of \textit{\bibinfo{series}{Graduate Texts in Mathematics}}.
\newblock \bibinfo{edition}{2nd} ed., \bibinfo{publisher}{Springer}, \bibinfo{address}{New York, NY}.
\newblock \DOIprefix\doi{10.1007/978-1-4612-0949-2}.
\bibitem[{Kessler et~al.(2012)Kessler, Lindner and Sorensen}]{kessler2012statistical}
\bibinfo{author}{Kessler, M.}, \bibinfo{author}{Lindner, A.}, \bibinfo{author}{Sorensen, M.}, \bibinfo{year}{2012}.
\newblock \bibinfo{title}{Statistical methods for stochastic differential equations}.
\newblock \bibinfo{publisher}{Chapman and Hall/CRC}.
\bibitem[{Lau and Siu(2008)}]{lau2008option}
\bibinfo{author}{Lau, J.W.}, \bibinfo{author}{Siu, T.K.}, \bibinfo{year}{2008}.
\newblock \bibinfo{title}{On option pricing under a completely random measure via a generalized esscher transform}.
\newblock \bibinfo{journal}{Insurance: Mathematics and Economics} \bibinfo{volume}{43}, \bibinfo{pages}{99--107}.
\bibitem[{Mandelbrot and Van~Ness(1968)}]{mandelbrot1968fractional}
\bibinfo{author}{Mandelbrot, B.B.}, \bibinfo{author}{Van~Ness, J.W.}, \bibinfo{year}{1968}.
\newblock \bibinfo{title}{Fractional brownian motions, fractional noises and applications}.
\newblock \bibinfo{journal}{SIAM review} \bibinfo{volume}{10}, \bibinfo{pages}{422--437}.
\bibitem[{McDonald and Xu(1995)}]{mcdonald1995generalization}
\bibinfo{author}{McDonald, J.B.}, \bibinfo{author}{Xu, Y.J.}, \bibinfo{year}{1995}.
\newblock \bibinfo{title}{A generalization of the beta distribution with applications}.
\newblock \bibinfo{journal}{Journal of Econometrics} \bibinfo{volume}{66}, \bibinfo{pages}{133--152}.
\bibitem[{Milevsky and Posner(1998)}]{milevsky1998asian}
\bibinfo{author}{Milevsky, M.A.}, \bibinfo{author}{Posner, S.E.}, \bibinfo{year}{1998}.
\newblock \bibinfo{title}{Asian options, the sum of lognormals, and the reciprocal gamma distribution}.
\newblock \bibinfo{journal}{Journal of financial and quantitative analysis} \bibinfo{volume}{33}, \bibinfo{pages}{409--422}.
\bibitem[{Monfort and Pegoraro(2012)}]{monfort2012asset}
\bibinfo{author}{Monfort, A.}, \bibinfo{author}{Pegoraro, F.}, \bibinfo{year}{2012}.
\newblock \bibinfo{title}{Asset pricing with second-order esscher transforms}.
\newblock \bibinfo{journal}{Journal of Banking \& Finance} \bibinfo{volume}{36}, \bibinfo{pages}{1678--1687}.
\bibitem[{Rombouts and Stentoft(2014)}]{rombouts2014bayesian}
\bibinfo{author}{Rombouts, J.V.}, \bibinfo{author}{Stentoft, L.}, \bibinfo{year}{2014}.
\newblock \bibinfo{title}{Bayesian option pricing using mixed normal heteroskedasticity models}.
\newblock \bibinfo{journal}{Computational Statistics \& Data Analysis} \bibinfo{volume}{76}, \bibinfo{pages}{588--605}.
\bibitem[{Schoutens(2003)}]{schoutens2003Levy}
\bibinfo{author}{Schoutens, W.}, \bibinfo{year}{2003}.
\newblock \bibinfo{title}{L{\'e}vy processes in finance: pricing financial derivatives}.
\newblock \bibinfo{publisher}{Wiley Online Library}.
\bibitem[{Stroock(2011)}]{stroock2011probability}
\bibinfo{author}{Stroock, D.W.}, \bibinfo{year}{2011}.
\newblock \bibinfo{title}{Probability Theory: An Analytic View}.
\newblock \bibinfo{edition}{2nd} ed., \bibinfo{publisher}{Cambridge University Press}, \bibinfo{address}{Cambridge}.
\bibitem[{Yu and Zhang(2005)}]{yu2005three}
\bibinfo{author}{Yu, K.}, \bibinfo{author}{Zhang, J.}, \bibinfo{year}{2005}.
\newblock \bibinfo{title}{A three-parameter asymmetric laplace distribution and its extension}.
\newblock \bibinfo{journal}{Communications in Statistics—Theory and Methods} \bibinfo{volume}{34}, \bibinfo{pages}{1867--1879}.
\bibitem[{Øksendal(2003)}]{Oksendal2003}
\bibinfo{author}{Øksendal, B.}, \bibinfo{year}{2003}.
\newblock \bibinfo{title}{Stochastic Differential Equations: An Introduction with Applications}.
\newblock \bibinfo{edition}{6th} ed., \bibinfo{publisher}{Springer}, \bibinfo{address}{Berlin, Heidelberg}.
\newblock \DOIprefix\doi{10.1007/978-3-642-14394-6}.

\end{thebibliography}

\end{document}